\documentclass[12pt]{article}

\usepackage[a4paper, total={6in, 8in}]{geometry}

\usepackage{amsmath,amssymb,amsthm}
\usepackage{stmaryrd}
\usepackage{dsfont}
\usepackage[utf8]{inputenc}
\usepackage[backend=biber, style=numeric-comp]{biblatex}
\usepackage{graphicx} 
\usepackage[hidelinks]{hyperref}
\usepackage{fullpage}
\usepackage{cancel}
\usepackage{enumerate}
\usepackage{todonotes}
\usepackage{soul}
\usepackage{float}
\usepackage{xcolor}
\usepackage{dsfont}
\usepackage{tikz}

\DeclareMathOperator*{\argmax}{argmax}
\newtheorem{result}{Result}

\usepackage{csquotes}

\newcommand{\densite}{{\mathrm{f}}}

\newcommand\plan{p}
\newcommand\Hamming{H}

\newcommand{\Prob}{\mathrm{P}}
\newcommand{\XU}{\mathcal{X}}
\newcommand{\obs}{\mathrm{obs}}

\newcommand{\Pop}{\mathcal{U}}
\usepackage{amsmath, amssymb}
\usepackage{graphicx} 
\usepackage{fancyhdr} 
\usepackage{hyperref} 

\title{Differential Privacy and Survey Sampling}
\author{Daniel Bonnéry, Julien Jamme}
\date{\today}

\begin{document}

\maketitle

\footnotetext[1]{Daniel, Bernard Bonnéry, daniel.bonnery@ensae.fr}
\footnotetext[2]{Julien Jamme, julien.jamme@insee.fr}

\maketitle

\begin{abstract}
    
The Horvitz-Thompson estimate of a total can be seen as as differentially private mechanism applied to this population total. We provide forumlae to compute the $\epsilon$ and $\delta$ parameter for this specific mecanism, coupled or not coupled with the addition of a Laplace or a Gaussian noise.
This allows to determine the scale of the Laplace privacy mechanism to be added to reach a specified level of privacy, expressed in terms of $\epsilon,\delta$ differential privacy.
In particular, we provide simple formulae for the special case of simple random sampling on binary data.

\end{abstract}

\section*{Introduction}
National Statistical Institutes invest significant resources in conducting surveys to achieve a desired level of accuracy in their estimates. However, if the $\epsilon,\delta$ differential pricacy  criterion were to be applied, these institutes would be required to introduce noise into their results, thereby compromising the very accuracy they have worked to secure. Moreover, an survey estimate can be seen as the true total to which a sampling noise has been added, that is the outcome of a privacy mechanism.
The question is how private is this mechanism ? We answer this question.
In section 1, we provide the statistical framework and introduce the key concepts.
In section 2, we provide our main results.
In section 3, we discuss the importance of the results and their application in real life situations.
The appendix contains the proofs.

\section{Definitions}

\subsection{Conventions and notations}
The random variables introduced in this article are defined on the same probability space, whose probability is denoted by $\Prob$, and take values in real vector spaces equipped with the Borel $\sigma$-algebra. For a random variable $X$, 
$\Prob^X$ denotes the image measure of  $\Prob$ by $X$. Capital letters are generally used for random variables, caligraphic capital letters for sets, lower case for non random variables.

\subsection{Population and Individual Characteristics}
Consider a set of statistical units ${\Pop}=\{i: i\in {\Pop}\}$ referred to as the  ``  population'' , of size $N$. For each statistical unit $i$, we denote by $X_{i}$ its one dimensional characteristic, taking values in a closed bounded set $\XU_1$ de $\mathbb{R}$,  whose minimum and maximum bounds are  $m_x$ and $M_x$. The vector of characteristics for all statistical units $X$ takes values in the set $\XU:=\left\{x\in\XU_1^\Pop:t(x) \in[m_t,M_t] \right\}$, where $t:x\mapsto \sum_{i\in {\Pop}}x_i$, and $N~m_x\leq m_t<M_t\leq N~M_x\in\mathbb{R}$.

\subsection{ Sampling and Survey Data}

We are given ${\plan}$ a sampling plan without replacement on $\Pop$, that is, a probability law on the measurable space formed by $\Pop$ and the $\sigma$-algebra generated by its subsets. The sample $S$ is a random variable with law  ${\plan}$, satisfying the property :
$\forall x\in\XU, P^{S\mid X=x}={\plan}$. We denote  ${\plan}_i={\plan}\left(\{i\in S\}\right)$, ${\plan}_{i,j}=\left(\{\{i,j\}\subset S\}\right)$,
${\plan}_{i,j,\ell}={\plan}\left(\{\{i,j,\ell\}\subset S\}\right)$, ${\plan}_{-i}=1-{\plan}_i$ (\cite Gourieroux). 
To a given sample  $s\subset\Pop$, we associate the survey data ($\mathrm{obs}_X(s)$), where for $x\in \XU$, $\mathrm{obs}_x(s)=\left(s,(x_i)_{i\in s}\right)$. The model $\{ {\plan}^{\mathrm{obs}_x}: x\in\XU\}$ is commonly called the fixed population model for plan-based inference. Under this model, the so-called Horvitz-Thompson estimator of the total  $t(x)$, is defined by:
$$\hat{t}(X,S)=\sum_{i\in S}\frac{X_i}{{\plan}_i}.$$

\subsection{Differential Privacy at threshold $\epsilon$, and at thresholds  ($\epsilon,\delta$)}

As a reminder, the Hamming distance on $\XU$ is defined by $\XU\times \XU\to \mathbb{R}, x,x'\mapsto \Hamming(x,x')=\sum_{i=1}^N (x_{i}\neq x'_{i})$. Let $x,x'\in\mathcal{\Hamming}_1(\XU)$, $i\in\Pop$ such that $x_i\neq x_i'$.
We denote $\mathcal{\Hamming}_1 :=\{(x,x')\in \XU: \Hamming(x,x')\leq 1\}$,  $\mathcal{\Hamming}_1^t:=\left\{\left(t(x),t(x')\right): (x,x')\in\mathcal{\Hamming}_1(\XU)\right\}$ and $\bar B_1(x):=\{x'\in\XU,\Hamming(x')\leq1\}$. We denote by $x_{-i}$ the vector $(x_j)_{j\in \Pop,j\neq i}$ and $(x_{-i},x'_{i})$ denotes the vector  $(x''_i)_{i\in \Pop}$ such that  $x''_{-i}=x_{-i}$ et $x''_i=x'_i$ .

A random variable $Z$ satisfies the differential privacy criterion at thresholds  $(\epsilon,\delta)$ for data from $\XU$ (\cite{dwork2006calibrating}), if and only if the condition :
$$\mathrm{DP}_{Z,\XU}(\epsilon,\delta ): =\left[\forall (x,x')\in\mathcal{\Hamming}_1(\XU) : DP_{Z,\XU}(\epsilon,\delta,x,x')\right],$$ is satisfied, with  $DP_{Z,\XU}(\epsilon,\delta,x,x'):=\left [\forall A\in\mathfrak{S}(Z),~P(A\mid X=x)\leq \exp(\epsilon) \times P(A\mid X=x')+\delta\right],$
and $\mathfrak{S}(Z)$ the $\sigma$-algebra induced by  $Z$.
We define: :
$\delta_{Z,\XU}(\epsilon):=\inf\{\delta\in\mathbb{R}^+: DP_{Z,\XU}(\epsilon,\delta)\},$ and
$\delta_{Z,\XU}(\epsilon,x,x'):=\inf\{\delta\in\mathbb{R}^+: DP_{Z,\XU}(\epsilon,\delta,x,x')\}$. 

We have: $\delta_{Z,\XU}(\epsilon)=\sup\left\{ \delta_{Z,\XU}(\epsilon,x,x'): (x, x')\in\mathcal{\Hamming}_1(\XU)\right\}.$ We are interested, in terms of differential privacy properties, in the following random variable : $Z_b=\hat{t}(X,S)+{bW}$, where ${W}$  is a random variable independent of $X$, and of  $S$, distributed according to a Laplace distribution with location parameter $0$and scale $1$. The law of $Z_b\mid X=x$ admits a density, denoted $\densite_{Z\mid X=x}$, with respect to the measure $\nu_b$, which is either the Lebesgue measure (pour $b>0$), or the counting measure on  $\mathbb{R}$ (for $b=0$). Under these conditions,
 $$\delta_{Z_b,\XU}(\epsilon,x,x'):=\int_{z\in\mathbb{R}}\left(\densite_{Z\mid X=x}(z)-\exp(\epsilon)\times \densite_{Z\mid X=x'}(z)\right)_+ \mathrm{d}\nu_b(z),$$ with  $(a)_+=a\times \mathds{1}_{\mathbb{R}^+}(a)$.
 We also have: $\delta_{Z_b,\XU}(\epsilon,x,x'):=P(Z_b\in A_{\epsilon,x,x'}\mid X=x)-e^\epsilon P(Z_b\in A_{\epsilon,x,x'}\mid X=x')$, with $A_{\epsilon,x,x'}=\{z:\densite_{Z_b\mid X=x}(z)-\exp(\epsilon)\times \densite_{Z_b\mid X=x'}(z)>0\}$.
 
 We define $\epsilon_{Z_b,\XU}(\delta,x,x')=\inf(\epsilon\in\mathbb{R}:\delta_{Z_b,\XU}(\epsilon,x,x')<\delta\}$, $\epsilon_{Z_b,\XU}(\delta)=\inf(\epsilon\in\mathbb{R}:\delta_{Z_b,\XU}(\epsilon)<\delta\}$, et $b_{\XU}(\delta,\epsilon):=\inf\left\{b\in[0,+\infty]: DP_{Z_b,\XU}(\delta,\epsilon)\right\} $.
  The objective of this article is to find the exact value or a lower bound of $b$.
 
\subsection{Conditional Differential Privacy given the Sample}

\citeauthor{drechsler2024complexitiesdifferentialprivacysurvey} (\cite{drechsler2024complexitiesdifferentialprivacysurvey}) iscuss the application of differential privacy methods in the context of survey data from surveys. Like Lin et al. (2023),\citeauthor{lin2023differentially}(\cite{lin2023differentially}),  the selection process is not viewed as a perturbative process.

The criterion sought is for a given realization $s$ de $S$:
$\forall A\in\mathfrak{S}(Z_b)$, $\forall y,y'\in \{\obs_x(s):x\in\XU\}$ tel que $H(y,y')\leq 1$, 
$$P(A\mid X_s=y,S=s)\leq \exp(\epsilon) \times P(A\mid X_s=y',S=s)+\delta.$$

This criterion is consistent with the need to protect a dataset (survey data) based on the principle that a potential attacker, to use the term used in the literature on data protection, knows the list of individuals selected for the survey. A less extreme scenario is considered here, where the attacker does not observe the sample $S$ (the list of selected statistical units).

\section{Results}

\subsection{In the absence of perturbation ($b=0$)}

For the first result, it should be noted that the data of $(X_{i})_{i\in S}$ alone is different from the data of the values $(X_{i})_{i\in S}$ and the index  $(i)_{i\in S}$.  The random vector  $(X_{i})_{i\in S}$ contains only unlabeled values by a population identifier $i$, as in a pseudonymized file.
The first result is interpreted as follows: with probability   $(1-\plan_i)$ the statistical unit $i$ does not appear in the sample and in this case no information is given on this unit. 

\begin{result}
$\delta_{\mathrm{obs}_X(S),\XU}(0)\leq\max_{i\in \Pop}\plan_i$.\label{resultat:dp1}
\end{result}

\begin{proof}
See Section \ref{sec:demores1}.
\end{proof}

We deduce from the previous result that any statistic from the survey data is $D.P.\left(0,\max_{i\in {\Pop}} {\plan}_i\right).$
We therefore have: $\delta(0)\leq\max_{i\in \Pop}\plan_i$.
The following result states that under certain conditions on  $\mathcal{X}$, $\delta(0)=\max_{i\in \Pop}\plan_i$,
which is interpreted as follows: without perturbation, for at least one value of  $x$, with a probability of $\max_{i\in \Pop} p_i$, 
which corresponds to the probability that the individual with maximum inclusion probability is in the sample, 
 knowing  $x_{-i}$ and $Z$, it is possible to know with certainty that $x_i$ belongs to a finite set that is a strict subset of $\XU(1)$.

\begin{result}[The inequality of  \ref{resultat:dp1}  is tight]\label{resultat:dp2}
    If $\XU_1$ is infinite, then $\delta_{\mathrm{obs}(X,S),\XU}(0)=\max\{{\plan}_j:j\in {\Pop}\}.$
\end{result}

If ${\plan}(\emptyset)<1$, and if $\mathrm{cardinal}(\mathcal{X}_1)>2^N$,  
then the privacy mecanism  associated with the scheme ${\plan}$ and the Horvitz-Thompson estimator is not  $D.P.\left(\epsilon,0\right)$.

\begin{result}[$\epsilon_{\hat{t}(X,S),\XU}(0)=+\infty$]\label{resultat:dp3}
If $\exists ~i\in {\Pop}, x\in \XU$, tel que  {$\mathrm{cardinal}\left(\left\{x'\in\XU:x'_{-i}=x_{-i}\right\}\right)>\mathrm{cardinal}(\{s\subset {\Pop} : {\plan}(s)>0, i\in s\}),$} then
$\epsilon_{\hat{t}(X,S),\XU}(0)=+\infty.$ 
\end{result}

The previous result is related to the fact that under certain conditions on  $\XU_1$ only, 
it is possible, with probability $p_i$,  to determine a finite list of possible values of $x_i$ knowing $x_{-i}$. However, this requires being able to calculate the values of the Horvitz-Thompson estimator for all samples with non-zero probability, and the number of such samples can be astronomical, as can the list of possible values.

\subsection{General case, addition of a Laplacian noise $(b\geq 0)$}

The random variable $Z_b$, conditionally on $X=x$ et $S=s$,  admits a density with respect to the measure  $\nu_b$ given by  $\densite_{Z_b\mid X=x,S=s}(z)= \densite_{bW}\left(z-\hat{t}(x,s)\right)$. By applying the formula for the marginal density, we deduce the density of $Z_b$ conditionally on $X=x$:
$$\densite_{Z_b\mid X=x}(z)=\sum_{s\subset {\Pop}}{\plan}(s)\times \densite_{bW}\left(z-\hat{t}(x,s)\right),$$ where
$\densite_{{bW}}(w)=(2b)^{-1}\exp\left(-\left|w\right|/b\right)$ if $b>0$, $
\mathds{1}_{\{0\}}(w)$ otherwise.
If $b>0$, then $\densite_{Z_b\mid X=x}$ is continuous.
The computation of  $\densite_{Z_b\mid X=x}$ is of order $O\left(\mathrm{cardinal}\left(\left\{s\subset \Pop: p(s)>0\right\}\right)\right)$ set aside the computation of $p(s)$.
Define: $\mathrm{supp}(g):=\{z\in\mathbb{R}:g(z)>0\}$ le support d'une fonction réelle $g$ positive.
Let

$m_{\hat{t}}:=\min\left(\bigcup_{x\in\XU}\mathrm{supp}(\densite_{Z_0\mid X=x})\right)$ et 
$M_{\hat{t}}:=\max\left(\bigcup_{x\in\XU}\mathrm{supp}(\densite_{Z_0\mid X=x})\right)$.

\begin{eqnarray*}
     \delta_{Z_0,\XU}(\epsilon,x,x')&=&
\sum_{z\in\mathbb{R}}
    \left(\sum_{s\subset\Pop}
        p(s) ~
            \left(
                \mathds{1}_{\{z\}}(\hat{t}(x,s))-
                \exp(\epsilon)~\mathds{1}_{\{z\}}(\hat{t}(x',s))
            \right)
    \right)_+\\&=&
\sum_{z\in\mathrm{supp}(\densite_{Z_0\mid X=x,S\ni i})}
    \left( \densite_{Z_0\mid X=x}(z)  
    -\exp(\epsilon)~
    \densite_{Z_0\mid X=x'}(z)
    \right)_+
\end{eqnarray*}
The theoretical formulae for $\delta_{Z_b,\XU}$ can lead to untractable compuations, and are given below:

\begin{eqnarray*}
     \lefteqn{\delta_{Z_b,\XU}(\epsilon,x,x')}\\
     &=&(2b)^{-1}
\int_{z\in\mathbb{R}}
    \left(\sum_{s\subset\Pop}
        p(s) 
            \left(\exp(-|z-\hat{t}(x,s)|/b)-
                \exp(\epsilon)\exp(-|z-\hat{t}(x',s)|/b))
            \right)
    \right)_+\mathrm{d}z\\&=&
    (2b)^{-1}
        \left(\sum_{s\subset\Pop}
                p(s)\left( 
                e^{-\hat{t}(x,s)/b}-
                e^{\epsilon -\hat{t}(x',s)/b}\right)
        \right)_+
    \int_{z=-\infty}^{m_{\hat{t}}}
  \exp(z/b)~\mathrm{d}z\\
        &&+
    (2b)^{-1}
        \left(\sum_{s\subset\Pop}
                p(s)\left( 
                e^{\hat{t}(x,s)/b}-
                e^{\epsilon+\hat{t}(x',s)/b)}\right)
        \right)_+
    \int_{z=M_{\hat{t}}}^{+\infty}
  \exp(-z/b)~\mathrm{d}z\\
        &&+
    (2b)^{-1}
\int_{z=m_{\hat{t}}}^{M_{\hat{t}}}
    \left(\sum_{s\subset\Pop}
        p(s) 
            \left(e^{-|z-\hat{t}(x,s)|/b}-
                e^ {\epsilon -|z-\hat{t}(x',s)/b|}
            \right)
    \right)_+\mathrm{d}z    
    \\&=& (1/2)
        \left(\sum_{s\subset\Pop}
                p(s)\left( 
                e^{-\hat{t}(x,s)/b}-
                e^{\epsilon -\hat{t}(x',s)/b}\right)
        \right)_+
  e^{m_{\hat{t}}/b}\\
        &&+
    (1/2)
        \left(\sum_{s\subset\Pop}
                p(s)\left( 
                e^{\hat{t}(x,s)/b}-
                e^{\epsilon+\hat{t}(x',s)/b)}\right)
        \right)_+
   e^{-M_{\hat{t}}/b}\\
        &&+
    (2b)^{-1}
\int_{z=m_{\hat{t}}}^{M_{\hat{t}}}
    \left(\sum_{s\subset\Pop}
        p(s) 
            \left(e^{-|z-\hat{t}(x,s)|/b}-
                e^ {\epsilon -|z-\hat{t}(x',s)/b|}
            \right)
    \right)_+\mathrm{d}z       
\\&=& (1/2)
        \left(\sum_{y\in\mathrm{supp}(\densite_{Z_0\mid X=x})}
                e^{-y/b}\left( \densite_{Z_0\mid X=x})\left(y\right)-
                e^{\epsilon}\densite_{Z_0\mid X=x'})\left(y\right)\right)                
        \right)_+
  e^{m_{\hat{t}}/b}\\
        &&+
    (1/2)
        \left(\sum_{y\in\mathrm{supp}(\densite_{Z_0\mid X=x})}
                e^{y/b}\left( \densite_{Z_0\mid X=x})\left(y\right)-
                e^{\epsilon}\densite_{Z_0\mid X=x'})\left(y\right)\right)              
        \right)_+
   e^{-M_{\hat{t}}/b}\\
        &&+
    (2b)^{-1}
\int_{z=m_{\hat{t}}}^{M_{\hat{t}}}
    \left(\sum_{y\in\mathrm{supp}(\densite_{Z_0\mid X=x})}
                e^{-|z-y|/b}\left( \densite_{Z_0\mid X=x})\left(y\right)-
                e^{\epsilon}\densite_{Z_0\mid X=x'})\left(y\right)\right)
    \right)_+\mathrm{d}z    
\end{eqnarray*}

\subsubsection{Determination of $\delta$ when $\epsilon=0$ is required and no additional noise is applied ($b=0$)}
This is a limit case, with theoretical interest, for which the expression of $\delta_{Z_0,\XU}(\epsilon,x,x')$ is simple:

\begin{result}{Expression of $\delta_{Z_0,\XU}(0,x,x')$}
\begin{eqnarray*}
    \delta_{Z_0,\XU}(0,x,x')&=&
\sum_{z\in\mathrm{supp}(\densite_{Z_0\mid X=x})}
    \left( \densite_{Z_0\mid X=x}(z)  
    -    \densite_{Z_0\mid X=x'}(z)
    \right)_+\\&=&
\left(\sum_{s\subset U, i\in s}
    p(s) \left(\mathds{1}\{z\}(\hat{t}(x,s))
    -    \mathds{1}\{z\}(\hat{t}(x',s))\right)
    \right)_+\\&\leq& p_i
\end{eqnarray*}
\end{result}

As a corollary,  if 
$\mathrm{supp}(\densite_{Z_0\mid X=x'})\cap\mathrm{supp}(\densite_{Z_0\mid X=x'})=\emptyset$ then $\delta_{Z_0,\XU}(0)=p_i$.

\subsubsection{Determination of $\epsilon_{Z_b,\XU}(0,x,x')$ when a Laplace noise is applied, and $\delta=0$ is required.}

Because $\epsilon_{Z_b,\XU}(\delta)=\sup\{\epsilon_{Z_b,\XU}(\delta,x,x'):(x,x)'\in\Hamming_1(\XU)$,
it is important to compute $\epsilon_{Z_b,\XU}(\delta,x,x')$.
We state a series of results that allow to compute  $\epsilon_{Z_b,\XU}(\delta,x,x')$ in specific cases.

\begin{result}
$$\epsilon_{Z_b,\XU}(0,x,x')=
    \left|
        \begin{array}{l}
            +\infty
                \text{ if }b=0\text{ et } \mathrm{supp}(\densite_{Z_0\mid X=x})\setminus\mathrm{supp}(\densite_{Z_0\mid X=x'})\neq \emptyset,\\
            \sup\left\{\densite_{Z_b\mid X=x}(z)/\densite_{Z_b\mid X=x'}(z):z\in\mathbb{R}\right\}  
                \text{ otherwise. }
        \end{array}
    \right.$$
\end{result}
\begin{proof}
Direct application of the definition of $\epsilon_{Z_b,\XU}(\delta,x,x')$.
\end{proof}

Define  $\hat{\mathcal{T}}_{x,x'}:=\{\hat{t}(x'',s):s\subseteq U,\plan(s)>0,x''\in\{x,x'\}\}$.

\begin{result}
If $b>0$ or $\mathrm{supp}(\densite_{Z_0\mid X=x})\subseteq
\mathrm{supp}(\densite_{Z_0\mid X=x'})$, then 
$$\sup\left\{\densite_{Z_b\mid X=x}(z)/
            \densite_{Z_b\mid X=x'}(z):z\in\mathbb{R}\right\}=\max\left\{\densite_{Z_b\mid X=x}(z)
            \densite_{Z_b\mid X=x'}(z):z\in \hat{\mathcal{T}}_{x,x'}\right\}.$$

\end{result}

\begin{proof}
Suppose $b>0$ or $\mathrm{supp}(\densite_{Z_0\mid X=x})\subseteq
\mathrm{supp}(\densite_{Z_0\mid X=x'})$. 

The ratio ${\densite_{Z_b\mid X=x}
            }/{
            \densite_{Z_b\mid X=x'}
            }$ restricted to each of the intervals  {$\left]-\infty, \min\left(\hat{\mathcal{T}}_{x,x'}\right)\right]$} and \hfill \break {$\left[\max\left(\bigcup_{s\subset \Pop, s\ni i}\{\hat{t}(x,s),\hat{t}(x,s')\}\right),+\infty\right[$} is constant.
            Indeed, 
for $z\in[-\infty,\min(\hat{\mathcal{T}}_{x,x'})] $,
if $b>0$ then :
 \begin{eqnarray*}
     \lefteqn{
        \frac{\densite_{Z_b\mid X=x}
            }{
            \densite_{Z_b\mid X=x'}
            }(z_{x,x'}(z))}\\&=&
            \frac{\sum_{s\subseteq \Pop }\plan(s) e^{z/b-\hat{t}(x,s)/b}}
            {\sum_{s\subseteq \Pop }\plan(s) e^{z/b-\hat{t}(x,s)/b}}\\&=&
            \frac{\sum_{s\subseteq \Pop }\plan(s) e^{-\hat{t}(x,s)/b}}
                 {\sum_{s\subseteq \Pop }\plan(s) e^{-\hat{t}(x',s)/b}}.
                 \end{eqnarray*}

For $z\in[\max(\hat{\mathcal{T}}_{x,x'}),+\infty] $,
then if $b>0:$
 \begin{eqnarray*}
     \lefteqn{
        \frac{\densite_{Z_b\mid X=x}
            }{
           \densite_{Z_b\mid X=x'}
            }(z_{x,x'}(z))}\\&=&
            \frac{\sum_{s\subseteq \Pop }\plan(s) e^{\hat{t}(x,s)/b-z/b}}
            {\sum_{s\subseteq \Pop }\plan(s) e^{\hat{t}(x,s)/b-z/b}}\\&=&
            \frac{\sum_{s\subseteq \Pop }\plan(s) e^{\hat{t}(x,s)/b}}
                 {\sum_{s\subseteq \Pop }\plan(s) e^{\hat{t}(x',s)/b}}.
                 \end{eqnarray*}
 Moreovere the ratio is continuous, as a ratio of sums of continuous strictly positive functions.
Thus :
\begin{eqnarray*}
\epsilon_{0,x,x'}&=&\sup
    \left\{
        \frac{
            \densite_{Z_b\mid X=x}
            }{
            \densite_{Z_b\mid X=x'}
            }(z):z\in\mathbb{R}\right\}\\&=&
\max
    \left\{
        \frac{
            \densite_{Z_b\mid X=x}
            }{
            \densite_{Z_b\mid X=x'}
            }(z):z\in\left[\min(\hat{\mathcal{T}}_{x,x'}),\max(\hat{\mathcal{T}}_{x,x'})\right]\
    \right\}.\end{eqnarray*}
If $b=0$ and  $\mathrm{supp}(\densite_{Z_0\mid X=x})\subseteq
\mathrm{supp}(\densite_{Z_0\mid X=x'})$, then:
$$\epsilon_{Z_b,\XU}(0,x,x')=
\max
    \left\{
        \frac{
            \densite_{Z_b\mid X=x}
            }{
            \densite_{Z_b\mid X=x'}
            }:z\in\mathrm{supp}(\densite_{Z_0\mid X=x})
    \right\}$$

Let $\hat{\mathcal{T}}_{x,x'}=
         \left\{\hat{t}(x'',s):
            s\in\mathrm{supp}(p),
            x''\in\{x,x'\}\right\}$.
The set $\hat{\mathcal{T}}_{x,x'}$ is  finite and let $J_{x,x'}$ be its cardinality. Let $z_{x,x'}: \llbracket 0, J_{x,x'}+1\rrbracket\to  \hat{\mathcal{T}}_{x,x',i}$ the increasing index  of $\hat{\mathcal{T}}_{x,x'}$.
Remark that $z_{x,x'}=z_{x',x}$.
Let $u^+_{x,x',b}(j):=\sum_{j'>0, j'\leq j} \densite_{Z_0\mid X=x,S\ni i}(z_{x,x'}(j))\times\exp(z_{x,x'}(j)/b)$, 
 $u^-_{x,x',b}(j):=\sum_{j'> j,j\leq j} \densite_{Z_0\mid X=x,S\ni i}(z_{x,x'}(j))\times\exp(-z_{x,x'}(j)/b)$, 
 $j_{x,x'}(z):=j$ if $z\in [z_{x,x'}(j),z_{x,x'}(j+1)[$, $0$ if   
 $z<z_{x,x'}(1)$, and $J_{x,x'}+1$ if  $z>z_{x,x'}(J_{x,x'})$. Then 
 $\densite_{Z_b\mid X=x,S\ni i }(z)=
\exp(z/b) ~u^-_{x,x',b}(j_{x,x'}(z)) +\exp(-z/b)~ u^+_{x,x',b}(j_{x,x'}(z))
$. 
The set $\hat{\mathcal{T}}_{x,x'}$ corresponds to the points of non-differentiability of $\densite_{Z_b\mid X=x}/\densite_{Z_b\mid X=x'}$.
For $1\leq j<J$, on the interval $[z_{x,x'}(j),z_{x,x'}(j+1)]$, $j_{x,x'}$ is constant and the derivative of $\densite_{Z_b\mid X=x}/\densite_{Z_b\mid X=x'}$ is of the same sign (o.s.s) as the piecewise constant: $u^-_{x,x',b}(j)\times u^+_{x',x,b}(j)-u^-_{x',x,b}(j)\times u^+_{x,x',b}(j)$. Indeed 
on this interval, $\densite_{Z_b\mid X=x}/\densite_{Z_b\mid X=x'}$ is equal to 
$(u^-~e^{z/b}+u^+~e^{-z/b})/(u'^-~e^{z/b}+u'^+~e^{-z/b})$, 
with $u^-:=u^-_{x,x',b}(j)$, $u^+:=u^+_{x,x',b}(j)$, $u'^-:=u^-_{x',x,b}(j)$, $u'+:=u^+_{x',x,b}(j)$.
Its derivative is 
\begin{eqnarray*}
\lefteqn{\left(\mathrm{d}(\densite_{Z_b\mid X=x}/\densite_{Z_b\mid X=x'})/(\mathrm{d}(z)\right)(z)}\\
&=&\left(\exp(z/b) ~u'^-~+~\exp(-z/b)~ u'^+\right)^{-2}
\\&&\quad
\times
\left(\left(\exp(z/b) ~u^- ~-~\exp(-z/b)~ u^+\right)
\left(\exp(z/b) ~u'^-~+~\exp(-z/b)~ u'^+\right)\right.\\
&&\quad
\quad\left.-\left(\exp(z/b) ~u^- ~+~\exp(-z/b)~ u^+\right)
\left(\exp(z/b) ~u'^-~ -~\exp(-z/b)~ u'^+\right)\right)\\
&=&2 \left(\exp(z/b) ~u'^-~+~\exp(-z/b)~ u'^+\right)^{-2}\times ~\left(u^- u'^+-u'^- u^+\right)\\
&o.s.s.&u^- u'^+-u'^- u^+.
\end{eqnarray*}
Therefore for $1\leq j<J$, $\densite_{Z_b\mid X=x}/\densite_{Z_b\mid X=x'}$ reaches its maximum on  $[z_{x,x'}(j),z_{x,x'}(j+1)]$ in $z_{x,x'}(j)$, or in $z_{x,x'}(j+1)$. Furthermore the ratio $\densite_{Z\mid X=x}/\densite_{Z\mid X=x'}$ is constant on  the interval
$]-\infty,z_{x,x'}(1)]$ on and also on the interval $]z_{x,x'}(J),+\infty,[$. Therefore if  $b>0$ or $\mathrm{supp}(\densite_{Z_0\mid X=x})\subseteq
\mathrm{supp}(\densite_{Z_0\mid X=x'})$, then
$$\argmax
    \left\{
        \frac{
            \densite_{Z_b\mid X=x}
            }{
            \densite_{Z_b\mid X=x'}
            }(z):z\in\mathbb{R}
    \right\}=
    \argmax
    \left\{
        \frac{
            \densite_{Z_b\mid X=x}
            }{
            \densite_{Z_b\mid X=x'}
            }(z):z\in\hat{\mathcal{T}}_{x,x'}
    \right\}.$$

\end{proof}

Define: 
\begin{eqnarray*}
    \mathcal{S}_{-i}^+&:=&\{s\subseteq\Pop,p(s)>0,S \not\ni i\}\\
    \mathcal{S}_{-i}^-&:=&\{s\subseteq\Pop,p(s)>0,S \not\ni i\}\\
    \mathcal{S}_{i}^{+}&:=&\{s\subseteq\Pop,p(s)>0,S \ni i,\hat{t}_{-i}(x,s)<z-x_i/\plan_i\}\\
    \mathcal{S}_{i}^{\mp}&:=&\{s\subseteq\Pop,p(s)>0,S \ni i,z-x'_i/\plan_i\leq\hat{t}_{-i}(x,s)<z-x_i/\plan_i\}\\
    \mathcal{S}_{i}^{-}&:=&\{s\subseteq\Pop,p(s)>0,S \ni i,z-x'_i/\plan_i\leq\hat{t}_{-i}(x,s)\}\\
\end{eqnarray*}

\begin{result}
If $\mathcal{S}_{-i}^-=\emptyset$ or $\mathcal{S}_{i}^+\emptyset$ , and $x_i\geq x'_i$ then
$$\argmax_{z\in\mathbb{R}}(\densite_{Z_b\mid X=x,S\ni i}/\densite_{Z_b\mid X=x',S\ni i})(z)=
\max\{\hat{t}_{-i}(x,s): s\subseteq\Pop,\plan(s)>0, s\ni i\}.$$
If $\mathcal{S}_{-i}^-=\emptyset$ or $\mathcal{S}_{i}^+\emptyset$ , and $x'_i\geq x_i$ then
$$\argmax_{z\in\mathbb{R}}(\densite_{Z_b\mid X=x,S\ni i}/\densite_{Z_b\mid X=x',S\ni i})(z)=
\max\{\hat{t}_{-i}(x,s): s\subseteq\Pop,\plan(s)>0, s\ni i\}.$$
\end{result}

\begin{proof}

We define the arbitrary order on the set of sets:
$\mathcal{S}_{-i}^-\prec
    \mathcal{S}_{-i}^+\prec
    \mathcal{S}_{i}^{-}\prec
    \mathcal{S}_{i}^{\mp}\prec
    \mathcal{S}_{i}^{+}$ and $\mathcal{S}:=\{\mathcal{S}_{-i}^-,
    \mathcal{S}_{-i}^+,
    \mathcal{S}_{i}^{-},
    \mathcal{S}_{i}^{\mp},
    \mathcal{S}_{i}^{+}\}$.

Define 
\begin{eqnarray*}
a(x,x',s,s')&:=&e^{(\hat{t}(x,s)-|z-\hat{t}(x',s')|)/b)}\\
c(x,x',s,s')&:=&a(x,x',s,s')-a(x',x,s,s')\\
d(x,x',s,s')&:=&c(x,x',s,s')+c(x,x',s',s)
 \\&=&e^{(\hat{t}(x ,s)-|z-\hat{t}(x',s')|)/b)}-
               e^{(\hat{t}(x',s)-|z-\hat{t}(x,s')|)/b)}\\
               &&+e^{(\hat{t}(x ,s')-|z-\hat{t}(x',s)|)/b)}-
               e^{(\hat{t}(x',s')-|z-\hat{t}(x,s)|)/b)}.
\end{eqnarray*}

\begin{eqnarray*}
 \lefteqn{
 \frac{
    \densite_{Z_b\mid X=x}
       }{
    \densite_{Z_b\mid X=x'}
    }(z_J)- 
    \frac{
            \densite_{Z_b\mid X=x}
            }{
            \densite_{Z_b\mid X=x'}
            }(z)}\\
            &=&
 \frac{\sum_{A\in\mathcal{S}}\sum_{s\in A} p(s)e^{-z_J/b}e^{\hat{t}(x,s)/b}}
 {\sum_{A\in\mathcal{S}}\sum_{s\in A} p(s)e^{-z_J/b}e^{\hat{t}(x',s)/b}}-
            \frac{\sum_{A\in\mathcal{S}}\sum_{s\in A} p(s)e^{-|z-\hat{t}(x,s)|/b}}
 {\sum_{A\in\mathcal{S}}\sum_{s\in A} p(s)e^{-|z-\hat{t}(x',s)|/b}}
            \\&o.s.s.a.&
 \sum_{A,B\in\mathcal{S}}\sum_{s\in A,s'\in B} 
 p(s)p(s')c(x,x',s,s')
            \\&=&
 \frac12\sum_{A\in\mathcal{S}}\sum_{s,s'\in A} 
 p(s)p(s')\left(d(x,x',s,s')\right)+
  \sum_{A,B\in\mathcal{S}\atop A\prec B}\sum_{s\in A\atop s'\in B} 
 p(s)p(s')\left(d(x,x',s,s')\right)
\end{eqnarray*}

We compute:

For $s,s'\in \mathcal{S}_{-i}^-, \hat{t}(x,s)=\hat{t}(x',s)$ and $\hat{t}(x,s')=\hat{t}(x',s')$ so
\begin{eqnarray*}
c(x,x',s,s')&=&e^{(\hat{t}(x ,s)+z-\hat{t}(x',s'))/b)}-
               e^{(\hat{t}(x',s)+z-\hat{t}(x,s'))/b)}\\
               &=&e^{(\hat{t}(x ,s)+z-\hat{t}(x',s'))/b)}-
               e^{(\hat{t}(x,s)+z-\hat{t}(x',s'))/b)}\\
               &=&0.
\end{eqnarray*}

For $s,s'\in \mathcal{S}_{-i}^+, \hat{t}(x,s)=\hat{t}(x',s)$ and $\hat{t}(x,s')=\hat{t}(x',s')$ so $c(x,x',s,s')=0$

For $s,s'\in \mathcal{S}_{i}^-, \hat{t}(x,s)=\hat{t}_{-i}(x,s)+x_i/\plan_i$ and 
$\hat{t}(x',s)=\hat{t}_{-i}(x,s)+x'_i/\plan_i$ so
\begin{eqnarray*}
c(x,x',s,s')&=&e^{(\hat{t}(x ,s)+z-\hat{t}(x',s'))/b)}-
               e^{(\hat{t}(x',s)+z-\hat{t}(x,s'))/b)}\\
            &=&e^z\left(e^{(x_i-x'_i)/(\plan_i b)}e^{(\hat{t}_{-i}(x ,s)-\hat{t}_{-i}(x,s'))/b}-
               e^{(x'_i-x_i)/(\plan_i b)}e^{(\hat{t}_{-i}(x,s)-\hat{t}_{-i}(x,s'))/b}\right)\\
d(x,x',s,s')&=&e^z\left(e^{(x_i-x'_i)/(\plan_i b)}
                    \left(e^{(\hat{t}_{-i}(x ,s)-\hat{t}_{-i}(x,s'))/b}
                            +e^{(\hat{t}_{-i}(x ,s')-\hat{t}_{-i}(x,s))/b}\right)\right.
                            \\&&\quad\quad\left.-
               e^{(x'_i-x_i)/(\plan_i b)}
               \left(e^{(\hat{t}_{-i}(x',s)-\hat{t}_{-i}(x,s'))/b}
               +e^{(\hat{t}_{-i}(x,s')-\hat{t}_{-i}(x,s))/b}\right)\right)\\
               &=&4~\mathrm{ch}((\hat{t}_{-i}(x ,s)-\hat{t}_{-i}(x,s'))/b) ~e^z~\mathrm{sh}\left({(x_i-x'_i)/(\plan_i b)}\right)\\
               &\geq&0.
\end{eqnarray*}

For $s,s'\in \mathcal{S}_{i}^\mp$, 
\begin{eqnarray*}
c(x,x',s,s')&=&e^{(\hat{t}(x ,s)-z+\hat{t}(x',s'))/b)}-
               e^{(\hat{t}(x',s)+z-\hat{t}(x,s'))/b)}\\
               &=&e^{(x_i+x'_i)/(\plan_i b)}~e^{-z}~e^{(\hat{t}_{-i}(x ,s)+\hat{t}_{-i}(x,s'))/b)}\\
               &&-e^{(x'_i-x_i)/(\plan_i b)}~e^z~
               e^{(\hat{t}_{-i}(x,s)-\hat{t}_{-i}(x,s'))/b)}\\
               &=&2~e^{(\hat{t}_{-i}(x ,s)+x'_i/\plan_i)/b)}~
               \left(\mathrm{sh}\left((\hat{t}_{-i}(x,s')+x_i/\plan_i-z)/b\right)\right)\\&\geq&0.
\end{eqnarray*}
For $s,s'\in \mathcal{S}_{i}^+$, \\
\begin{eqnarray*}
c(x,x',s,s')&=&e^{(\hat{t}(x ,s)+\hat{t}(x',s')-z)/b)}-
               e^{(\hat{t}(x',s)+\hat{t}(x,s')-z)/b)}\\
d(x,x',s,s')&=&0.
\end{eqnarray*}
For $s\in \mathcal{S}_{-i}^-,s'\in\mathcal{S}_{-i}^+$, $\hat{t}(x,s)=\hat{t}(x',s)$, so:
\begin{eqnarray*}
c(x,x',s,s')&=&e^{(\hat{t}(x ,s)-z+\hat{t}(x',s'))/b)}-
               e^{(\hat{t}(x',s)-z+\hat{t}(x,s'))/b)}\\
            &=&e^{(\hat{t}(x ,s)-z+\hat{t}(x,s'))/b)}-
               e^{(\hat{t}(x,s)-z+\hat{t}(x,s'))/b)}\\
               &=&0,\\
c(x,x',s',s)&=&e^{(\hat{t}(x ,s')+z-\hat{t}(x',s))/b)}-
               e^{(\hat{t}(x',s')+z-\hat{t}(x,s))/b)}\\
               &=&e^z-e^z=0,\\
d(x,x',s,s')&=&0.               
\end{eqnarray*}
For $s\in \mathcal{S}_{-i}^-,s'\in\mathcal{S}_{i}^-$, $\hat{t}(x,s)=\hat{t}(x',s)$,
\begin{eqnarray*}
c(x,x',s,s')&=&e^{(\hat{t}(x ,s)+z-\hat{t}(x',s'))/b)}-
               e^{(\hat{t}(x',s)+z-\hat{t}(x,s'))/b)}\\
            &=&\left(e^{-x'_i/(\plan_i~b)}-e^{-x_i/(\plan_i~b)}\right)~e^{(\hat{t}(x,s)+z-\hat{t}_{-i}(x,s'))/b)}\\
            &\geq&0,\\
c(x,x',s',s)&=&e^{(\hat{t}(x ,s')+z-\hat{t}(x',s))/b)}-
               e^{(\hat{t}(x',s')+z-\hat{t}(x,s))/b)}\\
            &=&\left(e^{x_i/(\plan_i~b)}-e^{x'_i/(\plan_i~b)}\right)~e^{(\hat{t}_{-i}(x,s')+z-\hat{t}(x,s))/b)}\\
            &\geq&0.
\end{eqnarray*}
For $s\in \mathcal{S}_{-i}^-,s'\in\mathcal{S}_{i}^\mp$
\begin{eqnarray*}
d(x,x',s,s')&=&e^{(\hat{t}(x ,s)-z+\hat{t}(x',s'))/b)}-
               e^{(\hat{t}(x',s)+z-\hat{t}(x,s'))/b)}\\
               &&+e^{(\hat{t}(x ,s')+z-\hat{t}(x',s))/b)}-
               e^{(\hat{t}(x',s')+z-\hat{t}(x,s))/b)}\\
            &=&e^{(\hat{t}(x ,s)-z+x'_i/\plan_i+\hat{t}_{-i}(x,s'))/b)}-
               e^{(\hat{t}(x,s)+z-x_i/\plan_i-\hat{t}_{-i}(x,s'))/b)}\\
               &&+e^{(x_i/\plan_i+\hat{t}_{-i}(x ,s')+z-\hat{t}(x,s))/b)}-
               e^{(x'_i/\plan_i+\hat{t}_{-i}(x,s')+z-\hat{t}(x,s))/b)}\\
            &=&e^{(\hat{t}_{-i}(x,s')+x'_i/\plan_i)/ b}\left(
               2 ~\mathrm{sh}(\hat{t}(x ,s)-z)/b)\right)\\
               &&+2 ~e^z~\mathrm{sh}\left(x_i/\plan_i+\hat{t}_{-i}(x ,s')-\hat{t}(x,s))/b\right)\\
            &\geq&0.
\end{eqnarray*}
because the arguments passed to  the $\mathrm{sh}$ functions are positive.

For $s\in \mathcal{S}_{-i}^-,s'\in\mathcal{S}_{i}^+$,
\begin{eqnarray*}
d(x,x',s,s')&=&e^{(\hat{t}(x ,s)-z+\hat{t}(x',s'))/b)}-
               e^{(\hat{t}(x',s)-z+\hat{t}(x,s'))/b)}\\
               &&+e^{(\hat{t}(x ,s')+z-\hat{t}(x',s))/b)}-
               e^{(\hat{t}(x',s')+z-\hat{t}(x,s))/b)}\\
               &=&\left(e^{x'_i/(\plan_i~b)}-e^{x_i/(\plan_i~b)}\right)~e^{(\hat{t}(x ,s)-z+\hat{t}_{-i}(x,s'))/b)}\\
               &&+\left(e^{x_i/(\plan_i~b)}-e^{x'_i/(\plan_i~b)}\right)~e^{(\hat{t}_{-i}(x ,s')+z-\hat{t}(x,s))/b)}\\
               &=&2~e^{\hat{t}_{-i}(x ,s')/b}\left(e^{x_i/(\plan_i~b)}-e^{x'_i/(\plan_i~b)}\right)~
                  \mathrm{sh}\left((z-\hat{t}(x,s))/b)\right)
                  \\&\leq&0.
\end{eqnarray*}

For $s\in \mathcal{S}_{-i}^+,s'\in\mathcal{S}_{i}^-$,
\begin{eqnarray*}
d(x,x',s,s')&=&e^{(\hat{t}(x ,s)+z-\hat{t}(x',s'))/b)}-
               e^{(\hat{t}(x',s)+z-\hat{t}(x,s'))/b)}\\
               &&+e^{(\hat{t}(x ,s')-z+\hat{t}(x',s))/b)}-
               e^{(\hat{t}(x',s')-z+\hat{t}(x,s))/b)}\\
               &=&\left(e^{-x'_i/(\plan_i~b)}-e^{-x_i/(\plan_i~b)}\right)
               e^{(\hat{t}_{-i}(x ,s)+z-\hat{t}_{-i}(x,s'))/b)}\\
               &&+\left(e^{x_i/(\plan_i~b)}-e^{x'_i/(\plan_i~b)}\right)
               e^{(\hat{t}_{-i}(x ,s')-z+\hat{t}_{-i}(x,s))/b)}\\
               &\geq&0.
\end{eqnarray*}
For $s\in \mathcal{S}_{-i}^+,s'\in\mathcal{S}_{i}^\mp$,
\begin{eqnarray*}
d(x,x',s,s')&=&e^{(\hat{t}(x ,s)-z+\hat{t}(x',s'))/b)}-
               e^{(\hat{t}(x',s)+z-\hat{t}(x,s'))/b)}\\
               &&+e^{(\hat{t}(x ,s')-z+\hat{t}(x',s))/b)}-
               e^{(\hat{t}(x',s')-z+\hat{t}(x,s))/b)}\\
            &=&e^{x'_i/(\plan_i b)}e^{(\hat{t}(x,s)-z+\hat{t}_{-i}(x,s'))/b)}-
               e^{-x_i/(\plan_i b)}e^{(\hat{t}(x,s)+z-\hat{t}_{-i}(x,s'))/b)}\\
               &&+e^{x_i/(\plan_i b)}e^{(\hat{t}_{-i}(x ,s')-z+\hat{t}(x,s))/b)}-
               e^{x'_i/(\plan_i b)}e^{(\hat{t}_{-i}(x,s')-z+\hat{t}(x,s))/b)}\\
            &=&2 e^{\hat{t}(x,s)/b}\mathrm{sh}(\left((x_i/\plan_i +\hat{t}_{-i}(x ,s')-z)/b)\right)\\
            &\geq&0.
\end{eqnarray*}
For $s\in \mathcal{S}_{-i}^+,s'\in\mathcal{S}_{i}^+$
\begin{eqnarray*}
d(x,x',s,s')&=&e^{(\hat{t}(x ,s)-z+\hat{t}(x',s'))/b)}-
               e^{(\hat{t}(x',s)-z+\hat{t}(x,s'))/b)}\\
               &&+e^{(\hat{t}(x ,s')-z+\hat{t}(x',s))/b)}-
               e^{(\hat{t}(x',s')-z+\hat{t}(x,s))/b)}\\
            &=&0.
\end{eqnarray*}
For $s\in \mathcal{S}_{i}^-,s'\in\mathcal{S}_{i}^\mp$

\begin{eqnarray*}
d(x,x',s,s')&=&e^{(\hat{t}(x ,s)-z+\hat{t}(x',s'))/b)}-
               e^{(\hat{t}(x',s)+z-\hat{t}(x,s'))/b)}\\
               &&+e^{(\hat{t}(x ,s')+z-\hat{t}(x',s))/b)}-
               e^{(\hat{t}(x',s')+z-\hat{t}(x,s))/b)}\\
               &=&e^{(x_i+x'_i)/(\plan_i b)}e^{(\hat{t}_{-i}(x ,s)-z+\hat{t}_{-i}(x,s'))/b)}
               -e^{(-x_i+x'_i)/(\plan_i b)}e^{(\hat{t}_{-i}(x,s)+z-\hat{t}_{-i}(x,s'))/b)}\\
               &&+e^{(x_i-x'_i)/(\plan_i b)}e^{(\hat{t}_{-i}(x ,s')+z-\hat{t}_{-i}(x,s))/b)}
               -e^{(-x_i+x'_i)/(\plan_i b)}e^{(\hat{t}_{-i}(x,s')+z-\hat{t}_{-i}(x,s))/b)}\\
               &=&e^{(\hat{t}_{-i}(x ,s)+x'_i/\plan_i )/b}~\mathrm{sh}\left((\hat{t}_{-i}(x ,s')+x_i/\plan_i -z)/b\right)\\
               &&+\mathrm{sh}\left((x_i-x'_i)/(\plan_i b)\right)~e^{(\hat{t}_{-i}(x ,s')+z-\hat{t}_{-i}(x,s))/b)}\\
               &\geq&0.
\end{eqnarray*}
For $s\in \mathcal{S}_{i}^-,s'\in\mathcal{S}_{i}^+$
\begin{eqnarray*}
d(x,x',s,s')&=&e^{(\hat{t}(x ,s)-z+\hat{t}(x',s'))/b)}-
               e^{(\hat{t}(x',s)-z+\hat{t}(x,s'))/b)}\\
               &&+e^{(\hat{t}(x ,s')+z-\hat{t}(x',s))/b)}-
               e^{(\hat{t}(x',s')+z-\hat{t}(x,s))/b)}\\
               &=&e^{(x_i+x'_i)/(\plan_i b)}e^{(\hat{t}_{-i}(x ,s)-z+\hat{t}_{-i}(x,s'))/b)}-
               e^{(x_i+x'_i)/(\plan_i b)}e^{(\hat{t}_{-i}(x,s)-z+\hat{t}_{-i}(x,s'))/b)}\\
               &&+e^{(x_i-x'_i)/(\plan_i b)}e^{(\hat{t}_{-i}(x ,s')+z-\hat{t}_{-i}(x,s))/b)}-
               e^{(-x_i+x'_i)/(\plan_i b)}e^{(\hat{t}_{-i}(x,s')+z-\hat{t}_{-i}(x,s))/b)}\\
               &=& 2~ \mathrm{sh}\left((x_i-x'_i)/(\plan_i b)\right)e^{(\hat{t}_{-i}(x ,s')+z-\hat{t}_{-i}(x,s))/b)}
               \\&\geq&0.
\end{eqnarray*}
For $s\in \mathcal{S}_{i}^\mp,s'\in\mathcal{S}_{i}^+$
\begin{eqnarray*}
d(x,x',s,s')&=&e^{(\hat{t}(x ,s)-z+\hat{t}(x',s'))/b)}-
               e^{(\hat{t}(x',s)-z+\hat{t}(x,s'))/b)}\\
               &&+e^{(\hat{t}(x ,s')-z+\hat{t}(x',s))/b)}-
               e^{(\hat{t}(x',s')+z-\hat{t}(x,s))/b)}    \\
               &=&e^{(x_i+x'_i)/(\plan_i b)}e^{(\hat{t}_{-i}(x ,s')-z+\hat{t}_{-i}(x,s))/b)}-
               e^{(-x_i+x'_i)/(\plan_i b)}e^{(\hat{t}_{-i}(x,s')+z-\hat{t}_{-i}(x,s))/b)}\\
               &=&2~e^{(\hat{t}_{-i}(x,s')+x'_i/\plan_i)/b}\mathrm{sh}\left((\hat{t}(x,s)-z)/b)\right)   \\
               &\geq&0.
\end{eqnarray*}

\end{proof}

\subsection{Cas particulier: données binaires et sondage aléatoire simple}
The number of combinaisons of $n$ between $N$ is $N \choose n $, equal to $(n!(N-n)!)^{-1}(N!)$ if $0\leq n\leq N$, and $0$ otherwise. Let $n, m_t,M_t\in\{1,\cdots,N\}$, $0\leq m_t < M_t\leq N$.
Supppose  ${\plan}(s)= 1/{N \choose n }$ if $\mathrm{cardinal}(s)=n$, $0$ otherwise. Besides, we assume that  $\XU=\left\{x\in\{0,1\}^N: \sum_i x_i\in\llbracket m_t,M_t\rrbracket\right\}$.

We then have, for 
$x\in \XU$:
$$\densite_{Z_b\mid X=x}(z)=
    \sum_{y=\max\{0,n+t(x)-N\}}^{\min\{t(x),n\} }
    \densite_{bW}\left(z-\frac{N}ny\right) 
    \frac{{t(x)\choose y}{{N-t(x)}\choose {n-y}}}{{N\choose n}}
.$$

This densitiy function does depend on  $x$ only via $t(x)$, and its computation is of complexity order $O(N^2)$, and so is the computation of  $\delta_{Z_b,\XU}(\epsilon)$ and $\epsilon_{Z_b,\XU}(\delta)$.

\begin{result}[Expression de $\delta_{Z_0,\XU}(\epsilon)$]
For simple random sampling for binary data:
If $b=0$, alors
$$\delta_{Z_0,\XU}(\epsilon)=
\max\left\{
\sum_{z=0}^n
\left(
    \frac{{t_x\choose z}{{N-t_x}\choose {n-z}}-\exp(\epsilon){t_x'\choose z}{{N-t_x'}\choose {n-z}}}{{N\choose n}}
\right)_+:(t_x,t_x')\in\mathcal{\Hamming}_1^t(\mathcal{X})\right\}$$
\end{result}

\begin{proof}
    \begin{eqnarray*}
    \delta_{Z_0,\XU},x,x'(\epsilon)&=&
\int
\left(
\sum_{y=0}^n
\densite_{0W}\left(z-\frac{N}ny\right)\times
    \frac{{t_x\choose y}{{N-t_x}\choose {n-y}}-\exp(\epsilon){t_x'\choose y}{{N-t_x'}\choose {n-y}}}{{N\choose n}}
\right)_+\mathrm{d}\nu_0(z)\\&=&
\sum_{z\in\mathbb{R}}
\left(
\sum_{y=0}^n
\mathds{1}_{\{0\}}\left(Nz/n-y\right)\times
    \frac{{t_x\choose y}{{N-t_x}\choose {n-y}}-\exp(\epsilon){t_x'\choose y}{{N-t_x'}\choose {n-y}}}{{N\choose n}}
\right)_+\\&=&
\sum_{z'=0}^n
\left(
   \frac{{t_x\choose z'}{{N-t_x}\choose {n-z'}}-\exp(\epsilon){t_x'\choose z'}{{N-t_x'}\choose {n-z'}}}{{N\choose n}}
\right)_+
\end{eqnarray*}

\end{proof}

The support of  $\sum_{i\in S}X_i$ conditionnally  $\sum_{i\in U} X_i=t_x$ is
$\llbracket \max\{0,n+t_x-N\},\min\{n,t_x\}\rrbracket$:

If $t_x'=t_x+1$, then $\llbracket \max\{0,n+t_x-N\},\min\{n,t_x\}\rrbracket=\llbracket \max\{0,n+t_x-N+1\},\min\{n,t_x+1\}\rrbracket$
$\Leftrightarrow \max\{(N-1)-t_x,t_x\}\geq n$.
Si $\max\{(N-1)-t_x,t_x\}\geq n$, then for  $z\in\llbracket 1,n\rrbracket$:
$\frac{{t_x\choose z}{{N-t_x}\choose {n-z}}}
    {{t_x'\choose z}{{N-t_x'}\choose {n-z}}}=
    \frac{(t_x+1-z)(N-t_x)}{(t_x+1)(N-t_x-n+z)}
    $

$$\lim_{z\to \delta\times\infty}
    \frac{\densite_{\hat{t}(x,S)+{W}}}
         {\densite_{\hat{t}(x',S)+{W}}}(z)=
    \frac{
        \sum_{s\subset {\Pop},s\not\ni i}{\plan}(s)\exp\left(\frac{\hat{t}(x,s)}{\sigma}\right)+
        \exp\left(\frac{\delta(x'_i-x_i)}{\sigma {\plan}_i}\right)
        \sum_{s\subset {\Pop},s\ni i}{\plan}(s)\exp\left(\frac{\hat{t}(x,s)}{\sigma}\right)
    }
    {
        \sum_{s\subset {\Pop},s\not\ni i}{\plan}(s)\exp\left(\frac{\hat{t}(x,s)}{\sigma}\right)+    
        \sum_{s\subset {\Pop},s\ni i}{\plan}(s)\exp\left(\frac{\hat{t}(x,s)}{\sigma}\right)
    }
    $$

\begin{result}
    \begin{eqnarray*}
    \lefteqn{\epsilon_{Z_0,\XU}(0)
    =}\\&&\left|\begin{array}{ll}+\infty &\text{ if }\min\{m_t,N-M_t\} < n\\
    \mathrm{ln}\left(
    \max\left\{
        \frac{N-M_t+1}{N-M_t+1-n},\frac{m_t+1}{m_t+1-n}
    \right\}
    \right) &\text{ otherwise }\end{array}\right.
    \end{eqnarray*}
\end{result}

\subsection{Additinal Gaussian privacy mechanism}

\begin{result}{Conditional to selection or non selection espected value and variance} 

For  $i\in {\Pop}$, if ${\plan}_i\notin\{0,1\}$, then $\forall x\in\mathcal{X}^{\Pop}$, 
\begin{align*}
\mathrm{E}\left[\hat{t}(X,S)\mid i\in S,X=x\right]&=    t_{-i|i}(x)+\frac{x_{i}}{{\plan}_i},
&\mathrm{Var}\left[\hat{t}(X,S)\mid i\in S,X=x\right]&=\sigma^2_{-i|i}(x)\\
\mathrm{E}\left[\hat{t}(X,S)\mid S\not\ni i,X=x\right]&=    t_{-i|-i}(x)
&\mathrm{Var}\left[\hat{t}(X,S)\mid S\not\ni i,X=x\right]&=\sigma^2_{-i|-i}(x).
\end{align*}

with 

\begin{align*}
I_i&=\mathds{1}_S(i),
    &I_{-i}&=1-I_i,\\
{\plan}_{j\mid i}&=\mathrm{E}\left[I_j\mid I_i=1\right]=\frac{{\plan}_{i,j}}{{\plan}_i},
    &{\plan}_{j\mid -i}&=\mathrm{E}\left[I_j\mid  I_i=0\right]=\frac{{\plan}_j-{\plan}_{i,j}}{1-{\plan}_i},\\
t_{-i|i}(x)&=\sum_{j\in {\Pop}\setminus\{i\}}\frac{x_j}{{\plan}_j} {\plan}_{j|i},
&t_{-i|-i}(x)&=\sum_{j\in {\Pop}\setminus\{i\}}\frac{x_j}{{\plan}_j}{\plan}_{j|-i},\\
{\plan}_{j,\ell\mid i}&=\mathrm{E}\left[I_jI_\ell\mid  I_i=1\right]=\frac{{\plan}_{i,j,\ell}}{{\plan}_i},
    &{\plan}_{j,\ell\mid -i}&=\mathrm{E}\left[I_jI_\ell\mid  I_i=0\right]=\frac{{\plan}_{j,\ell}-{\plan}_{i,j,\ell}}{1-{\plan}_i},\\
\sigma^2_{-i|i}(x)&=   \sum_{j,\ell\in {\Pop}\setminus\{i\}}\frac{x_jx_\ell}{{\plan}_j{\plan}_\ell}  ({\plan}_{j,\ell|i}-{\plan}_{j|i}{\plan}_{\ell|i}),
&\sigma^2_{-i|-i}(x)&= \sum_{j,\ell\in {\Pop}\setminus\{i\}}\frac{x_jx_\ell}{{\plan}_j{\plan}_\ell}  ({\plan}_{j,\ell|-i}-{\plan}_{j|-i}{\plan}_{\ell|-i}).
\end{align*}
\end{result}

Remark that quantities indexed by $-i$ : $t_{-i|-i}(x)$, $t_{-i|i}(x)$, $\sigma^2_{-i|i}(x)$, $\sigma^2_{-i|i}(x)$,   do not depend on $x_{i}$.

\begin{proof}
Compute : 
\begin{align*}
\end{align*}

\begin{eqnarray*}
\mathrm{E}\left[\hat{t}(X,S)\mid i\in S, X=x\right]
        &=&\sum_{j\in {\Pop}}\frac{x_j}{{\plan}_j} \mathrm{E}\left[I_j\mid  I_i=1\right]\\
        &=&\frac{x_{i}}{{\plan}_i}+\sum_{j\in {\Pop}\setminus\{i\}}\frac{x_j}{{\plan}_j} {\plan}_{j|i}  \quad \text{ car }{\plan}_{i|i}=1\\
        &=&\frac{x_{i}}{{\plan}_i}+t_{-i|i}(x)
\end{eqnarray*}

\begin{eqnarray*}
\mathrm{E}\left[\hat{t}(X,S)\mid X=x,S\not\ni i,{\plan}={\plan}\right]
        &=&\sum_{j\in {\Pop}}\frac{x_j}{{\plan}_j} \frac{{\plan}_j-{\plan}_{i,j}}{1-{\plan}_i}\\
        &=&\sum_{j\in {\Pop}\setminus\{i\}}\frac{x_j}{{\plan}_j} \frac{{\plan}_j-{\plan}_{i,j}}{1-{\plan}_i}\\
        &=&     t_{-i|-i}(x)
\end{eqnarray*}

\begin{eqnarray*}
    \mathrm{Var}\left[\hat{t}(X,S)\mid i\in S, X=x\right]
        &=&    \mathrm{Var}\left[\left.\frac{x_{i}}{{\plan}_i}+\sum_{j\in {\Pop}\setminus\{i\}}\frac{x_j}{{\plan}_j}\mathds{1}_S(j)\right| i\in S, X=x\right]\\
        &=&\sum_{j,\ell\in {\Pop}\setminus\{i\}}\frac{x_jx_\ell}{{\plan}_j{\plan}_\ell}\mathrm{Cov}\left[I_j,I_\ell\mid I_i=1,X=x,{\plan}={\plan}\right] \\
        &=&\sum_{j,\ell\in {\Pop}\setminus\{i\}}\frac{x_jx_\ell}{{\plan}_j{\plan}_\ell} ({\plan}_{j,\ell|i}-{\plan}_{j|i}{\plan}_{\ell|i})\\
        &=&\sigma^2_{-i|i}(x)
\end{eqnarray*}

\begin{eqnarray*}
\mathrm{Var}\left[\hat{t}_{-i}\mid S\not\ni i, X=x\right](x)
        &=&\sum_{j,\ell\in {\Pop}}\frac{x_jx_\ell}{{\plan}_j{\plan}_\ell}\mathrm{Cov}\left[I_j,I_\ell\mid I_i=0,X=x,{\plan}={\plan}\right] \\
        &=&\sum_{j,\ell\in {\Pop}\setminus\{i\}}\frac{x_jx_\ell}{{\plan}_j{\plan}_\ell} ({\plan}_{j,\ell|-i}-{\plan}_{j|-i}{\plan}_{\ell|-i})\\
        &=&    \sigma^2_{-i|-i}(x)
\end{eqnarray*}

\end{proof}

\section{Discussion}

Our results allow to compute the level of privacy of survey statistics.
It also provides the ability to use sampling as a privacy mechanism.

\appendix

\end{document}